\documentstyle[12pt]{article}    % Basic style
\newskip\stdskip                      % standard vertical space
\newcommand{\stdspace}{\hskip 0.75em plus 0.15em \ignorespaces}
\newcommand{\ppar}{\par\goodbreak\vskip 8pt plus 3pt minus 3pt} 

\newcommand{\rk}[1]{\ppar{\bf #1}\stdspace}    
\newenvironment{prf}{\par{\bf Proof}\ }{\hfill$\Box$\par}
\newcommand{\cl}{\centerline}         %  Centerline
\input xy
\xyoption{all}

\newtheorem{thrm}{Theorem}%[section]
\newtheorem{lem}[thrm]{Lemma}
\newtheorem{cor}[thrm]{Corollary}
\newtheorem{rem}[thrm]{Remark}
\newtheorem{defn}[thrm]{Definition}
\newtheorem{exmpl}[thrm]{Example}
\newtheorem{prop}[thrm]{Proposition}

\def\maprt#1{\, \smash{\mathop{\longrightarrow}\limits^{#1}}\, }
\def\maplft#1{\smash{\mathop{\longleftarrow}\limits^{#1}}}

\def\bd{  \begin{diagram}    }
\def\ed{  \end{diagram}      }

\def\thm #1  {\medskip\noindent{\bf #1}\quad}
\def\pf #1 {\smallskip\noindent{\bf #1}\par\nobreak\noindent}
\def\term #1{{\bf #1}}

\DeclareSymbolFont{AMSb}{U}{msb}{m}{n}
\DeclareMathSymbol{\ZZ}{\mathbin}{AMSb}{"5A}
\def\s{\Sigma}
\def\smsh{\wedge}
\def\om{\Omega}

\def\cross{\times}
\def\wdg{\vee}

\def\twdl{\widetilde}
\def\<{\langle}
\def\>{\rangle}

\def\sseq{\subseteq}

\def\bprp{\begin{prop}}
\def\eprp{\end{prop}}

\def\bthm{\begin{thrm}}
\def\ethm{\end{thrm}}

\def\blem{\begin{lem}}
\def\elem{\end{lem}}

\def\bcor{\begin{cor}}
\def\ecor{\end{cor}}

\def\brmk{\begin{rem}}
\def\ermk{\end{rem}}

\def\bdfn{\begin{defn}}
\def\edfn{\end{defn}}

\def\bexm{\begin{exmpl}}
\def\eexm{\end{exmpl}}

\def\|{\, \bigm|\, }

\def\A{{\cal A}}

\def\F{{\cal F}}
\def\M{{\cal M}}
\def\R{{\cal R}}
\def\S{{\cal S}}

\def\cl{{\mathrm cl}}

\def\holim{{\mathrm holim}}

\def\map{{\mathrm map}}

\def\!{}

\begin{document}
\setlength{\abovedisplayskip}{\stdskip}     % Reduces the space
\setlength{\belowdisplayskip}{\stdskip}     % around displays.
%
%
%                          Title page  
%                          ==========
%
%      Acknowledgements should not appear on this page. 
%      Please place these at the end of your introduction.
%
%
% This page will be reformatted by the Journal.  It will speed up the
% publication of your article if you type the title page information
% onto the standard form below.
%
% The following macros will turn this information into a basic format for 
% the page.  You are welcome to adapt these macros in any way you like.
%
%    Define the various ingredients of the title page:
%
\def\title#1{\def\thetitle{#1}}
\def\authors#1{\def\theauthors{#1}}
\def\author#1{\def\theauthors{#1}}
\def\address#1{\def\theaddress{#1}}
\def\secondaddress#1{\def\thesecondaddress{#1}}
\def\email#1{\def\theemail{#1}}
\def\url#1{\def\theurl{#1}}
\long\def\abstract#1\endabstract{\long\def\theabstract{#1}}
\def\primaryclass#1{\def\theprimaryclass{#1}}
\def\secondaryclass#1{\def\thesecondaryclass{#1}}
\def\keywords#1{\def\thekeywords{#1}}
%
%    Knuth's \ifundefined macro (needed to check for optional items):
\def\ifundefined#1{\expandafter\ifx\csname#1\endcsname\relax}
%
%   Basic title page layout (edit this macro if you
%   wish to adjust the title page layout) :
%
\long\def\maketitlepage{    % start of definition of \maketitlepage

\vglue 0.2truein   % top margin

% title :
%
{\parskip=0pt\leftskip 0pt plus 1fil\def\\{\par\smallskip}{\Large
\bf\thetitle}\par\medskip}   

\vglue 0.15truein 

% authors :
%
{\parskip=0pt\leftskip 0pt plus 1fil\def\\{\par}{\sc\theauthors}
\par\medskip}%
 
\vglue 0.1truein 

% address(es) email's and URL's (with switches to detect whether the
% optional items have been used) :
%
{\parskip=0pt\small
{\leftskip 0pt plus 1fil\def\\{\par}{\sl\theaddress}\par}
\ifundefined{thesecondaddress}\else\cl{and}  % second address?
{\leftskip 0pt plus 1fil\def\\{\par}{\sl\thesecondaddress}\par}\fi
\ifundefined{theemail}\else  % email address?
\vglue 5pt \def\\{\ \ {\rm and}\ \ } 
\cl{Email:\ \ \tt\theemail}\fi
\ifundefined{theurl}\else    % URL given?
\vglue 5pt \def\\{\ \ {\rm and}\ \ } 
\cl{URL:\ \ \tt\theurl}\fi\par}

\vglue 7pt 

{\bf Abstract}

\vglue 5pt

\theabstract

\vglue 7pt 

{\bf AMS Classification numbers}\quad Primary:\quad \theprimaryclass

Secondary:\quad \thesecondaryclass

\vglue 5pt 

{\bf Keywords:}\quad \thekeywords

%\np  % page break at the end of the title page

}    % end of definition of \maketitlepage
%
%
%   End of macros for basic title page layout
%
%
%  The following lines are for journal use.  Please do not disturb them.
%
%\input gtoutput
%\volumenumber{}\papernumber{}\volumeyear{}
%\pagenumbers{}{}\published{}
%\shorttitle{}  %\shortauthors{}  % for headlines (if needed)
%\proposed{}\seconded{}
%\received{}%\revised{}
%\accepted{}
%
%                  Title page information
%                  ======================
%
%   Type your title page information on the form below following
%   the format of the example.
%
%    \\ is the standard separator (between lines in \title and
%    \address, between authors and email addresses or URL's).
%
%
% Example:  \title{A short spoof paper\\with a two-line title}
% =======   \authors{Albert Einstein\\Leonardo da Vinci}
%           \address{IAS, Princeton}\secondaddress{Renaissance\\Venice}
%           \email{ae@ias.princeton.edu\\ldv@ren.ven.hist}
%           \abstract 
%           A short spoof paper with a very short abstract.
%           \endabstract 
%           \primaryclass{00-01, 00-02}\secondaryclass{68-00, 68-01}
%           \keywords{Short, spoof, paper}
%
%
%                  Start of title page form
%                  ========================
%
%    Type title, author(s) and address between the curly brackets:

\title{Miller Spaces and Spherical Resolvability of Finite Complexes} 

\authors{Jeffrey Strom} 
                 
\address{Dartmouth College, Hanover, NH 03755\\ 
{\tt{Jeffrey.Strom@Dartmouth.edu}} \\
{\tt{www.math.dartmouth.edu/\~{}strom/}} }              

%  second address, email address and URL (web address), are
%  all optional, uncomment if needed :

%\secondaddress{  }             
%\email{Jeffrey.Strom@Dartmouth.edu}                     
%\url{www.math.dartmouth.edu/\~{}strom/}                       

\abstract   % type your abstract below
We show that if $K$ is a nilpotent finite complex,
then $\om K$ can be built from spheres using fibrations and 
homotopy (inverse) limits.  This is applied to show
that if $\map_*(X,S^n)$ is weakly contractible
for all $n$, then $\map_*(X,K)$ is weakly contractible for
any nilpotent finite complex $K$.

\endabstract

%  AMS classification numbers, primary and secondary, and keywords :

\primaryclass{55Q05}    %Homotopy groups, general              
\secondaryclass{55P50}  %Localization and completion            

\maketitlepage

%
%%%%%%%%%%%%%%%%%%%%   End of title page
%
%%%%%%%%%%%%%%%%%%%%   Start of main body of article

\section*{Discussion of Results}

A  {\bf Miller space} is a CW complex $X$ with the property that 
the space of pointed maps from $X$ to $K$ 
is weakly contractible for every nilpotent
finite complex $K$, written $\map_*(X,K)\sim *$.  
They are named for Haynes Miller, who
proved in \cite{Miller} that the spaces 
$B\kern -.25em \ZZ\kern -.25em /p$ are
all Miller spaces; in fact, he proved that 
$\map_*(B\kern -.19em \ZZ\kern -.25em /p, K)$ 
is weakly contractible for
every finite dimensional CW complex $K$.  

In the stable category, one can define a {\it Miller
spectrum} by requiring that the mapping spectrum
$F(X,K)$ is contractible for every finite spectum $K$. 
Since cofibrations and fibrations are
the same in the stable category, a finite spectrum $K$ with $m$
cells is the fiber in a fibration
$K \maprt{} L \maprt{} S^n$ in which $L$ has only $m-1$ cells;
in the terminology of \cite{Cohen,Levi}, this means that
$K$ is {\it spherically resolvable with weight $m$}.
An easy induction shows that $X$ is a Miller spectrum if
and only if $F(X,S^n)\simeq *$ for every $n$.  

Our goal is to prove the following unstable analog of this 
observation:  if $\map_*(X,S^n)\sim *$, for all $n$,
then $X$ is a Miller space.  The proof of the stable version is not
available to us  because cofibrations are not fibrations,
unstably. To prove our result, it is necessary
to determine the extent to which a finite complex can
be constructed from spheres in a more general way, i.e., 
by arbitrary homotopy (inverse) limits 
\cite{Bousfield-Kan} and extensions by fibrations.   

To be more precise, we require some new
terminology.  We call a nonempty class $\R$
of spaces a \term{resolving class} if it is closed under 
weak equivalences and pointed homotopy (inverse) limits 
(all spaces, maps and homotopy limits will be pointed). It is a 
\term{strong resolving class} if it is further closed  
under extensions by fibrations, i.e., if whenever 
$F\maprt{} E\maprt{}B$ is a fibration with $F,B \in\R$, 
then $E\in\R$.
Resolving classes are dual to closed classes as defined in
\cite{Chacholski} and \cite[p.\thinspace 45]{EDF}.  

Notice that every resolving class $\R$ contains the one-point
space $*$ (cf. \cite[p.\thinspace 47]{EDF}).
From this, it follows that  if $F\maprt{}E\maprt{}B$
is a fibration with $E,B\in\R$, then $F\in\R$.  Similarly,
if $A_\alpha\in\R$ for each $\alpha$ then the
{\it categorical  product} $\Pi_\alpha A_\alpha\in\R$ also. 
The {\it weak product} $\twdl \Pi_\alpha A_\alpha$ 
is the homotopy colimit of the finite subproducts; if
for each $i$ only 
finitely many of the groups $\pi_i(A_\alpha)$ are nonzero, 
then the weak product has the same weak homotopy type as 
the categorical product. 

Let $\S$ be the smallest resolving class that contains $S^n$ for
each $n$, and let $\overline \S$ be the smallest strong resolving
class that contains $S^n$ for each $n$.  We say that a space
$K$ is \term{spherically resolvable} if $\om^k K\in \overline \S$
for some $k$.  This concept is related to, but not the same as, 
the notion  of spherical resolvability described in \cite{Cohen,Levi}.

\rk{Examples}     
\begin{enumerate}

\item[(a)]  If $f:A\maprt{}B$ is any map then the class of all $f$-local 
spaces is a resolving class \cite[p.\thinspace 5]{EDF}.  This includes, 
for example, the class of all spaces with $\pi_i(X) = 0$ for $i> n$, 
or all $h_*$-local spaces, where $h_*$ is a homology theory.  

\item[(b)]  If $P$ is a set of primes, then the class of all
$P$-local spaces is a strong resolving class. 

\item[(c)]  If $f: W\maprt{}*$, then the class of all
$f$-local spaces is a strong resolving class \cite[p.\thinspace 5]{EDF}.
This includes, for example, the class  $\{ K^+ \}$,
where $K^+$ denotes the Quillen plus construction on $K$ 
\cite[p.\thinspace 27]{EDF}.

\item[(d)] More generally, if $F$ is a covariant functor that
commutes with homotopy limits (and hence with fibrations) and $\R$ is a (strong)
resolving class, then the class 
$
\{ K \, |\, F(K) \in \R  \}
$
is also a (strong) resolving class.   This applies, for example
to the functor $F(K) = \map_*(X,K)$.

\item[(e)]  The class $\{ K \, | \, K\sim * \}$ is a strong resolving class.

\end{enumerate}

Our proofs will proceed by induction on a certain kind of cone length
\cite{A-S-S}.
Let $\F$ denote the collection of all finite type wedges of 
spheres.  The {\bf $\F$-cone length}  $\cl_\F(K)$ of a space $K$ is
the least integer $n$ for which there are cofibrations 
$S_i\maprt{} K_i \maprt{} K_{i+1}$,
$0\leq i < n$, with $K_0\simeq *$, $K_n\simeq K$
and each $S_i\in \F$.  If no
such $n$ exists, then $\cl_\F(K) =\infty$.  
Clearly every finite complex $K$ has $\cl_\F(K) < \infty$.

We denote by $\s\F\sseq \F$ the subcollection of all simply-connected
finite type wedges of spheres.
Finally, let $\S^\wdg$ be the smallest strong resolving class
that contains $\s\F$.

With these preliminaries in place, we can state our main
result.

\bthm\label{thrm:sres}
If $K$ is a nilpotent space with $\cl_\F(K) =n <\infty $, then
\begin{enumerate}
\item[{\rm (a)}]  $ K \in  \S^\wdg $,   
\item[{\rm (b)}]  $\om K\in \overline \S $, and
\item[{\rm (c)}]  $\om^n K \in  \S $.
\end{enumerate}
In particular, {\rm (b)} implies that 
every nilpotent finite complex $K$ is spherically
resolvable in our sense.
\ethm

Our application to Miller spaces follows from the following
more general consequence of Theorem \ref{thrm:sres}.
 
\bthm\label{thrm:main}
Let $\R$ be a strong resolving class and let $F$ be a functor
that commutes with homotopy limits. 
\begin{enumerate}
\item[{\rm (a)}]  Assume that $F(S^n)\in \R$ for each $n$.
Then $F(\om K)\in R$ for each nilpotent space $K$ with $\cl_\F(K) < \infty$.

\item[{\rm (b)}]  Assume that $F(S)\in \R$ for each $S\in\s\F$.
Then $F(K) \in \R$ for each nilpotent space $K$ with $\cl_\F(K) < \infty$.
\end{enumerate}
\ethm

To apply part  (b), we require the following result of Dwyer \cite{Dwyer}.

\bprp\label{prop:SinF}
Let $F$ be a functor that commutes with homotopy limits,
let $W$ be a space
and let $\R = \{ K\, |\, \map_*(W,F(K)) \sim *\}$. 
If $S^n\in\R$ for each $n$, then $\s\F\sseq \R$.
\eprp   

Together, Theorem \ref{thrm:main}(b) and Proposition \ref{prop:SinF}
immediately imply the desired statement about  Miller spaces.

\bcor\label{cor:appl}  If   $\map_*(X,S^n)\sim *$ for all $n$, 
then $\map_*( X,K)\sim *$ for every nilpotent space
$K$ with $\cl_\F(K) < \infty$.  In other words, $X$ is a Miller space.
\ecor

Corollary \ref{cor:appl} is by no means the only corollary of 
interest.  Other consequences are easily obtained by applying
Theorem \ref{thrm:main} to various strong resolving classes.  For example,
if $\map_*(X,S^n)$ is $P$-local for all $n$, then 
$\map_*(\s X,K)$ is $P$-local for every nilpotent space
$K$ with $\cl_\F(K) < \infty$.

If $X$ is simply-connected then Corollary \ref{cor:appl} can be strengthened
somewhat.  If $L$ is a space with a nilpotent covering space $K$
having $\cl_\F(K)<\infty$, then it is easy to see that $\map_*(X,L)\sim *$.
We end by making the surprising observation that a 
(non-nilpotent, of course) finite complex can be 
a Miller space!   

\noindent{\bf Example}\ \  
Let $A$ be a connected $2$-dimensional acyclic finite complex.
(The classifying space of the Higman group \cite{Hig} is such a space \cite{D-V};
so is the space obtained by removing a point from a homology $3$-sphere). 
Since $\pi_1(A)$ is equal to its commutator subgroup, there are no 
nontrivial homomorphisms from $\pi_1(A)$ to any nilpotent group.
It follows that if $f:A\maprt{} K$ with $K$ a nilpotent finite complex,
then $\pi_1(f) =0$ and so $f$ factors through $q:A\maprt{}A/A_1\simeq \bigvee S^2$.
Since $[A,S^2] \cong H^2(A) =0$, we conclude $f\simeq *$.
Thus $A$ is a Miller space.

\medskip

This example shows that the nilpotency hypothesis on the targets 
in Corollary \ref{cor:appl} cannot be entirely removed.  
It remains possible, however, that if $X$ is simply-connected 
and $\map_*(X,S^n)\sim *$ for each $n$, then $\map_*(X,K)\sim *$
for {\it every} finite complex $K$, or even for every finite-dimensional
complex.

{\it Acknowledgements}\  \ 
I would like to thank Robert Bruner and Charles McGibbon for
suggesting that I think about Miller spaces.   This work 
owes much to McGibbon in particular -- Corollary \ref{cor:appl}
was conjectured in joint work with him.
Thanks to Bill Dwyer for directing me to the result
of \cite{Hopkins}, which is the key to 
Proposition \ref{prop:desuspend}, and for the statement
and proof of Proposition \ref{prop:SinF}; thanks are also due to Daniel
Tanr\'e for bringing Proposition \ref{prop:hofiber}
to my attention.

\bigskip

\section{Proof of Theorem \ref{thrm:sres}}

We begin with two supporting results.

\bprp\label{prop:desuspend}
Let $K$ be a connected nilpotent space, let $\R$
be a resolving class and let $F$ be a functor that 
commutes with homotopy (inverse) limits.  
If  $F(\bigvee_{i=1}^m \s K)\in\R$
for each $m$, then $F(K)\in \R$.
\eprp

\begin{prf}
This follows from a result of Hopkins \cite[p.\thinspace 222]{Hopkins},
which says that $K$ is homotopy equivalent to
the homotopy (inverse) limit of a tower 
$$
A_0 \maplft{} A_1  \maplft{} \cdots \maplft{} A_n \maplft{} A_{n+1}\maplft{} \cdots
$$
of spaces, each of which is a homotopy (inverse) limit 
of a diagram of spaces of the form $\bigvee_{i=1}^m \s K$.
\end{prf}

%It follows that to show that $\om^k K \in\R$ it suffices to show 
%that $\om^k ( \bigvee_{i=1}^m\s K) \in \R$ for all $m$;
%more generally, if $F( \om^k ( \bigvee_{i=1}^m\s K))\in \R$
%for all $m$ then $F(\om^k K)\in\R$.

\bprp\label{prop:hofiber}
Let $A\maprt{}B\maprt{}C$ be a cofibration,
and let $F$ be the homotopy fiber of $B\maprt{}C$.
Then
$$
\s F \simeq \s A \wdg (\s A \smsh \om C).
$$
\eprp

\begin{prf}
Convert  the maps $A\maprt{*}C$, $B\maprt{}C$ and
$C\maprt{=}C$ to  fibrations.  The total spaces and fibers
form the commutative diagram
$$
\xymatrix{
A\cross \om C \ar[rd]\ar[rr]\ar[dd] && F \ar'[d][dd]\ar[rd]\\
& \om C \ar'[d][dd]\ar[rr] && {*}\ar[dd]\\
A \ar[rr]\ar[rd] && B \ar[rd]\\
& {*}  \ar[rr] && C,  \\
\\  }
$$ 

\vskip -.5in

\noindent
in which the bottom square is a homotopy pushout. 
A result of V. Puppe \cite{Puppe} shows that the top 
square is also a homotopy pushout.
Hence, the cofiber $\s F$ of the map
$F\maprt{} *$ has the same homotopy type as
the cofiber of $A\cross \om C\maprt{}\om C$, 
namely $\s A \wdg (\s A \smsh \om C)$, as can be seen from
the   diagram
$$
\xymatrix{
A\cross \om C \ar[d]\ar[r]\ar@{}[rd]|{{\mathrm pushout}}
& \om C \ar[r]\ar[d] & \s F\ar@{=}[d] \\
A \ar[r]^(.4){*} & A * \om C\ar[r] & \s A \wdg (\s A \smsh \om C) . }
$$ 
\end{prf}

\noindent{\bf Proof of Theorem \ref{thrm:sres}}\ \ 
Notice that the assumption on $X$ implies that $X$
is connected; we may therefore assume that $K$ is also
connected.

We prove  assertion (a) by induction on 
$\cl_\F(K)$.  If $\cl_\F(K)=1$, then $K$ is weakly equivalent
to a connected finite type wedge of spheres.
Therefore each $\bigvee_{i=1}^m \s  K\in\s\F$ 
and Proposition \ref{prop:desuspend} proves the assertion 
in the initial case.

Now assume that the result is known for all nilpotent 
spaces with $\F$-cone length less than $n$, and that $K$ is nilpotent
with $\cl_\F(K) = n$.  Two applications of Proposition \ref{prop:desuspend}  
reveal that it is enough to show $\bigvee_{i=1}^m \s^2 K  \in\S^\wdg$
for each $m$.

Write $V= \bigvee_{i=1}^m \s^2 K$.
Notice that $\cl_\F(\bigvee_{i=1}^m K) \leq \cl_\F(K)$,
and the double suspension of an $\F$-cone decomposition
of $\bigvee_{i=1}^m K$ is an $\F$-cone decomposition of $V$.
Thus we may assume that $V$ has an $\F$-cone decomposition
$S_i\maprt{} V_i\maprt{} V_{i+1}$, $0\leq i < n$    
with  $S_i, V_i\in \s\F$ for each $i$.  
Therefore,  we have a cofibration 
$L \maprt{} V \maprt{} W$ with $L$ simply-connected,
$\cl_\F(L) < n$ and $W\in\s\F$.
Let $F$ denote the homotopy fiber of $V\maprt{}W$,
so
$$
F \maprt{}V\maprt{}W
$$
is a fibration.  
Since $W \in \S^\wdg$, it suffices to show that $F\in \S^\wdg$.

Now we use Proposition \ref{prop:hofiber} to determine 
the homotopy type of $\s F$:
$$
\s F \simeq \s L \wdg ( L \smsh \s\om W) \simeq L\smsh 
\biggl(\bigvee_{\alpha} S^{n_\alpha}\biggr)
 $$
which is a finite type wedge of suspensions of $L$.   If we smash
an $\F$-cone length decomposition of $L$ with the space
$\bigvee_{\alpha} S^{n_\alpha}$ we obtain an $\F$-cone length
decomposition for $\s F$ -- in other words, $\cl_\F(\s F) < n$
and, more importantly, $\cl_\F(\bigvee_{i=1}^l \s F) < n$ for each $l$.

By the inductive hypothesis, $\bigvee_{i=1}^l \s F \in \S^\wdg$ 
for each $l$.   Since $L, V$
and $W$ are each simply-connected, so is $F$, and Proposition 
\ref{prop:desuspend} implies that $F\in \S^\wdg$, as desired.

To prove (b), observe that the collection $\M$ of
all $K$ with $\om K\in\overline \S$ is a resolving class
that contains $\s F$ by the Hilton-Milnor theorem \cite{G}.
Hence $\M$ contains all nilpotent spaces $K$ with $\cl_\F(K)<\infty$
by part (a). 

The proof of (c) is similar to the proof of (a).
The initial case of the induction is a special case of (b).  
To prove the inductive step, we write
$V= \bigvee_{i=1}^m \s^2 K$ and show that $V\in \S$.  As before, we 
consider the cofiber sequence $L\maprt{} V\maprt{}W$ with $W\in\s\F$ and the 
corresponding fibration $F\maprt{}V\maprt{}W$.  This gives us
a fibration
$$
\om^n V \maprt{} \om^n W \maprt{} \om^{n-1} F 
$$
with $\om^n W\in\S$. It now suffices to prove that $\om^{n-1} F\in\S$, 
which follows by induction using Proposition \ref{prop:desuspend}.
\hfill$\Box$\par

\medskip

\section{Proof of Theorem \ref{thrm:main} and Proposition \ref{prop:SinF}}

\noindent{\bf Proof of Theorem \ref{thrm:main}}\\ 
Let $\M$ be the class of all spaces $K$ such that 
$F(K) \in \R$; we have already seen
that $\M$ is a strong resolving class.  For part (a), 
$S^n\in \M$ for each $n$ by assumption, so $\overline \S \sseq \M$.
By Theorem \ref{thrm:sres}(b), $\M$ contains   $\om K$
for every nilpotent space  $K$ with $\cl_\F(K) < \infty$.
In part (b), we find that $\S^\wdg\sseq \M$, and so 
$\M$ contains every nilpotent space $K$ with $\cl_\F(K) < \infty$.
\hfill$\Box$\par

\medskip

\noindent{\bf Proof of Proposition \ref{prop:SinF}}\ \
Define a relation $<$ on $\s\F$ as follows:  $S<T$ if either
(1) the connectivity of $S$ is greater than the connectivity of $T$,
or
(2) the connectivity of $S$ equals the connectivity of $T$ (say
both are $(n-1)$-connected) and the rank of $\pi_n S$ is less than
the rank of $\pi_n T$. 
 
The key to this proof is the following claim.

\smallskip

\noindent {\sc Claim}\ \ Suppose that $S\in\F$   is
$(n-1)$-connected and has $\pi_n S\neq 0$. Then there is a map $f:S\maprt{} S^n$ 
such that the homotopy fibre $T$ of $f$ belongs to $\F$ and $T<S$. 

\noindent {\sc Proof of Claim}\ \ 
Write $S\sim S'\wdg S^n$, and let $f:S\maprt{} S^n$ be 
the map which collapses $S'$.  By \cite{G}, the homotopy fibre of $f$ is
$$
(S'\times\Omega S^n)/(*\times\Omega S^n)\sim S'\smsh(\Omega S^n)_+ \sim  
\bigvee_{m=0}^\infty \s^{(n-1)m} S' \in\s\F,
$$
using the James splitting of $\s\Omega S^n$.

\smallskip 

Now let $S=S_0 \in \s\F$.  Define $S_{n+1}$ as the 
fiber of a map $f:S_n\maprt{}S^{k(n)}$ as in the claim.  The result 
is a tower of spaces
$$
S_0\maplft{} S_1\maplft{} \cdots\maplft{} S_n\maplft{} S_{n+1}\maplft{}\cdots
$$
with $S_n\in\F$ and $S_{n+1} < S_n$ for each $n$.
Since spaces in this tower become arbitrarily highly connected
as $n$ increases, $\holim_n\, S_n \sim *$.

The fibrations $S_{n}\maprt{}S_{n+1}\maprt{} S^{k(n)}$ give rise
to fibrations
$$
\map_* (W, F(S_{n+1}))\maprt{}
\map_*(W, F(S_n)) \maprt{} 
\overbrace{\map_*(W, F(S^{k(n)}))}^{*}.
$$
It follows by induction that each map $S_n\maprt{} S$ induces a weak equivalence
$\map_*(W,F(S_n))\sim \map_*(W,F(S))$. Finally, we  compute
\[
\begin{array}{rcl}
\map_*(W,F(S))
&\sim&\holim_n\,\map_*(W,F(S))
\\
&\sim&\holim_n\,\map_*(W,F(S_n ))
\\
&\sim& \map_*(W,\holim_n\, F(S_n ))
\\
&\sim&\map_*(W,*)
\\
&\sim& {*}\\
\end{array}
\]
\hfill$\Box$\par

%%%%%%%%%%%%%%%%%%%%   End of main body of article
%
%                             References
%

%
%    Type any appendix material (to go after the references) here 
%
\end{document}